\documentstyle[12pt]{amsart}
\newcommand{\g}{\goth}

\newcommand{\gtsl}{\mbox{\g sl}}

\newcommand{\hgtsl}{\mbox{$\hat{\gtsl}$}}

\newcommand{\nz}{\mbox{${\Bbb Z}$}}

\newcommand{\vphi}{\varphi}
\theoremstyle{plain}
 \newtheorem{thm}{Theorem}
 \newtheorem{prop}[thm]{Proposition}
 \newtheorem{lemma}[thm]{Lemma}
 \newtheorem{cor}[thm]{Corollary}
\theoremstyle{definition}
 \newtheorem{defn}[thm]{Definition}

\theoremstyle{remark}

\begin{document}
\title{Quantum current operators (III)\\
Commutative quantum current operators, semi-infinite construction
and functional models\\} 
\author{Jintai Ding}
\author{Boris Feigin}
\address{Jintai Ding, RIMS, Kyoto University}
\maketitle
\begin{abstract}

We construct a commutative current operator 
$\bar x^+(z)$ inside  $U_q\left(\hat{\frak sl}(2)\right)$. 
With this operator and  the condition of quantum integrability 
on the quantum current of $U_q\left(\hat{\frak sl}(2)\right)$, we 
derive the quantization of the semi-infinite 
construction of integrable modules  of $\hat{\frak sl}(2)$ with 
the current operator $e(z)$ of   $\hat{\frak sl}(2)$. The 
quantization of the functional 
models for  $\hat{\frak sl}(2)$  are also given.

\end{abstract}
\pagestyle{plain}
\section{Introduction.} 

For any integrable highest weight module of $\hat{\frak sl}(2)$ of
level $k$, the current operators  $e(z)$ and $f(z)$ satisfy the following 
relations:
$$e(z)^{k+1}=f(z)^{k+1}=0,$$
which we call the condition of integrability \cite{LP}. 
For the case of quantum affine algebras, Drinfeld 
presented a formulation of affine quantum groups with generators 
in the form of current operators\cite{Dr2}, which, for the case of
  $U_q\left(\hat{\frak sl}(2)\right)$,   give us the 
 quantized current operators corresponding to $e(z)$ and $f(z)$ of 
 $\hat{\frak sl}(2)$. In \cite{DM}, we derive the quantum integrable condition
for  $U_q\left(\hat{\frak sl}(2)\right)$. 
On any level $k$ integrable module of 
$U_q\left(\hat{\frak sl}(2)\right)$ the matrix coefficients of 
$x^+(z_1)x^+(z_2)....x^+(z_{k+1})$ zero
 at $z_2/z_1=z_3/z_2=\ldots=z_{k+1}/z_k=q^2$,
and those of $x^-(z_1)x^-(z_2)....x^-(z_{k+1})$
are zero at $z_1/z_2=z_2/z_3=\ldots=z_k/z_{k+1}=q^2$,
where $x^+(z)$ and $x^-(z)$ are the quantized 
current operators of $U_q\left(\hat{\frak sl}(2)\right)$
corresponding to $e(z)$ and $f(z)$ of $\hat{\frak sl}(2)$, respectively. 
In the case of $\hat{\frak sl}(2)$, the integrable condition was used 
by Feigin and Stoyanovsky \cite{FS} to construct a level $k$ module
as a semi-infinite and 
describe the function models for the dual spaces. 
With the condition of the quantum integrability, 
still we can not simply  derive the quantization of the semi-infinite 
construction, because of the noncommutativity of the current operator
$x^+(z)$. Thus we have to modify the current operator $x^+(z)$ to 
 ``force'' it to commute  with itself. We use the subalgebra 
coming from the Heisenberg algebra of $U_q\left(\hat{\frak sl}(2)\right)$
to construct a commutative current operator 
$\bar x^+(z)=\Sigma \bar x_iz^{-i}$, which commutes with itself and 
$$  \bar x^+(z_1)\bar x^+(z_1q^2)....\bar x^+(z_{k+1}q^{2k})=0.$$
 Then the quantization of the  
semi-infinite construction, simply follows, namely 
the integrable modules of $U_q\left(\hat{\frak sl}(2)\right)$
can be identified with the space consisting of semi-infinite  expressions 
$\bar x^+_{i_1}...\bar x^+_{in}.....$, whose tails 
stablize in certain way and the action of $\bar x^+_i$ acts 
by multiplication.  Due to the introduction of 
the parameter $q$, we can  describe the action of  the operators
rigorously, especially the action of the operator 
$a_{-1}$, which corresponds to the operator $h_{-1}$ of 
$ \hat{\frak sl}(2)$. As in the case of \cite{FS}, 
the  functional models for the dual spaces of the subspace 
generated by $\bar x^+(z)$ on the highest weight vector 
in any irreducible integrable module of level $k$
are derived by using symmetric functions 
$f(t_1,....t_n)$, which is zero when $t_1=t_2q^2...=t_{k+1}q^{2k}.$

\section{$U_q\left(\hat{\frak sl}(n)\right)$ and commutative quantum 
current operators }
For the case of affine quantum  groups, Drinfeld gave a realization of
those algebras in terms of operators in the form of current\cite{Dr2}. 
We will first present such a realization for the case of 
$U_q\left(\hat{\frak sl}(n)\right)$.

Let $A=(a_{ij})$ be the Cartan matrix of type $A_{n-1}$.
\begin{defn}
The algebra $U_q(\hgtsl_n)$ is an associative algebra with unit 
1 and the generators: $\vphi_i(-m)$,$\psi_i(m)$, $x^{\pm}_i(l)$, for 
$i=i,...,n-1$, $l\in \nz $ and $m\in \nz_{\geq 0}$ and a central
 element $c$. Let $z$ be a formal variable and 
 $x_i^{\pm}(z)=\sum_{l\in \nz}x_i^{\pm}(l)z^{-l}$, 
$\vphi_i(z)=\sum_{m\in -\nz_{\geq 0}}\vphi_i(m)z^{-m}$ and 
$\psi_i(z)=\sum_{m\in \nz_{\geq 0}}\psi_i(m)z^{-m}$. In terms of the 
formal variables, 
the defining relations are 
\begin{align*}
& \vphi_i(z)\vphi_j(w)=\vphi_j(w)\vphi_i(z), \\
& \psi_i(z)\psi_j(w)=\psi_j(w)\psi_i(z), \\
& \vphi_i(z)\psi_j(w)\vphi_i(z)^{-1}\psi_j(w)^{-1}=
  \frac{g_{ij}(z/wq^{-c})}{g_{ij}(z/wq^{c})}, \\
& \vphi_i(z)x_j^{\pm}(w)\vphi_i(z)^{-1}=
  g_{ij}(z/wq^{\mp \frac{1}{2}c})^{\pm1}x_j^{\pm}(w), \\
& \psi_i(z)x_j^{\pm}(w)\psi_i(z)^{-1}=
  g_{ij}(w/zq^{\mp \frac{1}{2}c})^{\mp1}x_j^{\pm}(w), \\
& [x_i^+(z),x_j^-(w)]=\frac{\delta_{i,j}}{q-q^{-1}}
  \left\{ \delta(z/wq^{-c})\psi_i(wq^{\frac{1}{2}c})-
          \delta(z/wq^{c})\vphi_i(zq^{\frac{1}{2}c}) \right\}, \\
& (z-q^{\pm a_{ij}}w)x_i^{\pm}(z)x_j^{\pm}(w)=
  (q^{\pm a_{ij}}z-w)x_j^{\pm}(w)x_i^{\pm}(z), \\
& [x_i^{\pm}(z),x_j^{\pm}(w)]=0 \quad \text{ for $a_{ij}=0$}, \\
& x_i^{\pm}(z_1)x_i^{\pm}(z_2)x_j^{\pm}(w)-(q+q^{-1})x_i^{\pm}(z_1)
  x_j^{\pm}(w)x_i^{\pm}(z_2)+x_j^{\pm}(w)x_i^{\pm}(z_1)x_i^{\pm}(z_2) \\
& +\{ z_1\leftrightarrow z_2\}=0, \quad \text{for $a_{ij}=-1$}
\end{align*}
where
\[ \delta(z)=\sum_{k\in \nz}z^k, \quad
   g_{ij}(z)=\frac{q^{a_{ij}}z-1}{z-q^{a_{ij}}}\quad \text{about $z=0$} \]
\end{defn}
We define a grading on this algebra such that 
$x_i^\pm(n)$, $\vphi_i(n)$ and $\psi_i(n)$ are of degree $n$. 
We also always assume that $q$ is generic. 

Clearly, 
we have that $x^+(z)$ does not commute with itself. In order to do 
modify this operator, we have to rewrite the operators $\vphi_i(z)$ and 
$\psi_i(z)$ with new operators $a_{i,n}$ for $n\in \Bbb Z$. 

For the 
case of  $U_q(\hgtsl_2)$, we have that 
$$ \vphi(z)=
        \exp[-(q-q^{-1})\sum_{k\geq 0}a_{-k}z^k], $$
$$\psi(z)= 
        \exp[(q-q^{-1})\sum_{k\leq 0}a_{-k}z^{k}].$$

The new operators are defined as: 
$$ \exp[-(q-q^{-1})a_0]=  \vphi(0), $$ 
$$-(q-q^{-1})\sum_{k\geq 0}a_{-k}z^k= \log(1+  (\vphi(z)\vphi(0)^{-1}-1))= $$
$$\Sigma_{n> 0} (\vphi(z)\vphi(0)^{-1}-1)^n/n,  $$
$$(q-q^{-1})\sum_{k\leq 0}a_{k}z^{-k}= \log (1+(\psi(z)\psi(0)^{-1}-1))=$$
 $$\Sigma_{n>0} (\psi(z)\psi(0)^{-1}-1)^n/n.$$

\begin{prop}
$$[a_k,a_l]= \delta_{k+l, 0} (q^{2k}-q^{-2k})(q^c-q^{-c})/(k(q-q^{-1})^2),$$
$$[a_k, x^{\pm}(l)]=  (q^{2k}-q^{-2k})q^{\mp|c|/2}x^{\pm}(k+l)/(k(q-q^{-1})) . $$
\end{prop}

Let $k^-(z)$ be an current operator in $U_q\left(\hat{\frak sl}(2)\right)$,
 such that 
\begin{align*}
\text{(I)} & k^-(z)= 1+\Sigma_{n>0} k^-(n)z^{-n}, \\
\text{(II)} &k^-(z)x^+(w)=\frac {z-wq^2}{z-w} x^+(w)k^-(z), 
\end{align*}
where $k^-(n)$ are operators of degree $
n$, $k^-(n)k^-(m)=k^-(m)k^-(n)$ and 
 $\frac {z-wq^2}{z-w}$ is expanded in $w/z$. 

 Let $\bar x^+(w)= x^+(w)k^-(w)$, then we have 
\begin{prop}
$$(z-w) \bar x^+(z)\bar x^+(w)= (z-w)\bar x^+(w)\bar x^+(z). $$
\end{prop}

The proof comes from the following calculation. 
 $$(z-w) \bar x^+(z)\bar x^+(w)
= x^+(z) k^-(z)x^+(w)k^-(w)=(z-w)\frac {z-wq^2}{z-w}  x^+(z) x^+(w) k^-(z)k^-(w)
= $$
$$  ({zq^2-w})x^+(w) x^+(z)k^-(w) k^-(z)
= x^+(w)k^-(w) x^+(z) k^-(z)  ({zq^2-w}) (\frac {w-zq^2}{w-z})^{-1}=$$
$$ (z-w)\bar x^+(w)\bar x^+(z)$$

\begin{thm}
$$\bar x^+(z)\bar x^+(w)= \bar x^+(w)\bar x^+(z). $$
\end{thm}

{\bf Proof} 
Let $V_k$ be an integrable module of $U_q\left(\hat{\frak sl}(2)\right)$ 
and $V_k^*$ be its dual. Let $v\in V_k$ and $v^*\in V_k^*$. 
From the commutation relation between 
$x^+(z)$ and $x^+(w)$, we have that as a analytic function, 
the matrix coefficient $<v^*,  x^+(z)x^+(w)v>$ is zero when 
$z=w$, thus the correlation function $<v^*,  x^+(z)x^+(w)v>$ always has 
a factor $z-w$. This implies that the correlation function of 
 $<v^*,  \bar x^+(z)\bar x^+(w)v>$ does not have a pole at the 
point $z=w$, thus 
$$\bar x^+(z)\bar x^+(w)= \bar x^+(w)\bar x^+(z) $$ follows from 
$$(z-w)\bar x^+(z)\bar x^+(w)= (z-w)\bar x^+(w)\bar x^+(z). $$

\begin{prop}
Let
$$k^-(z)= \exp[(q-q^{-1}) \Sigma_{n>0}-q^{2(n+c/2}/(1+q^{2n}) a_n z^{-n}].$$
Then $ k^-(z)$ satisfies (I) and (II). 
\end{prop}

From now on, we will denote 
$x^+(z) \exp[(q-q^{-1}) \Sigma_{n>0}-q^{2(n+c/2}/(1+q^{2n}) a_n z^{-n}]$ by $\bar x(z)$
through out this paper. 

Let $k^+(z)$ be an current operator in $U_q\left(\hat{\frak sl}(2)\right)$
 such that 
\begin{align*}
\text{(I')} & k^+(z)= 1+\Sigma_{n<0} k^+(n)z^{-n},  \\
\text{(II')}& k^+(z)x^+(w)=\frac {z-w}{zq^2-w} x^+(w)k^+(z),
\end{align*}
where $k^+(n)$ are operators of degree $n$, $k^+(n)k^+(m)=k^+(m)k^+(n)$ and 
 $\frac {z-w}{zq^2-w}$ is expanded in $z/w$. 
 Let $\tilde x^+(w)=k^+(w)\tilde x^+(w)$, then we have 
\begin{prop}
$$(z-w) \tilde x^+(z)\tilde x^+(w)= (z-w)\tilde x^+(w)
\tilde x^+(z). $$
\end{prop}

The proof comes from the following calculation. 
 $$(z-w) \tilde x^+(z)\tilde x^+(w)
= k^+(z) x^+(z)k^+(w) x^+(w)=(z-w)\frac {z-wq^2}{z-w} k^+(z)k^+(w) 
 x^+(z) x^+(w) = $$
$$  ({zq^2-w})k^+(z)k^+(w) x^+(w) x^+(z)
=k^+(w) x^+(w) k^+(z)x^+(z) ({zq^2-w}) (\frac {z-w}{zq^2-w})=$$
$$ (z-w)\tilde x^+(w)\tilde x^+(z)$$

\begin{thm}
$$\tilde x^+(z)\tilde x^+(w)= \tilde x^+(w) \tilde x^+(z). $$
\end{thm}

The proof is the same as that if the theorem above. 

\begin{prop}
Let 
$$k^+(z)= \exp[-(q-q^{-1})\Sigma _{n<0}-q^{2(n+c/2}/(1+q^{2n}) a_nz^{-n}].
$$ Then 
$k^+(z)$ satisfies the condition (I') and (II'). 
\end{prop}. 

From now on, we will denote the operator $ \exp[-(q-q^{-1})\Sigma _{n<0}-q^{2(n+c/2}/(1+q^{2n}) a_nw^{-n}] x^+(w)$ by 
$\tilde x^+(w)$.

For the 
case of  $U_q(\hgtsl_n)$, we have that 
$$ \vphi_i(z)=
        \exp[-(q-q^{-1})\sum_{k\geq 0}a_{i,-k}z^k], $$
$$\psi_i(z)= 
        \exp[(q-q^{-1})\sum_{k\leq 0}a_{i,-k}z^{k}].$$

The new operators are defined as: 
$$ \exp[-(q-q^{-1})a_{i,0}]=  \vphi_i(0), $$ 
$$-(q-q^{-1})\sum_{k\geq 0}a_{i,-k}z^k= \log(1+  (\vphi_i(z)\vphi_i
(0)^{-1}-1))= $$
$$\Sigma_{n> 0} (\vphi_i(z)\vphi_i(0)^{-1}-1)^n/n,  $$
$$(q-q^{-1})\sum_{k\leq 0}a_{i,k}z^{-k}= \log (1+(\psi_i(z)\psi_i
(0)^{-1}-1))=$$
 $$\Sigma_{n>0} (\psi_i(z)\psi_i(0)^{-1}-1)^n/n.$$

Let 
$$k_i^+(z)= \exp[-(q-q^{-1})\Sigma _{n<0}-q^{2(n+c/2}/(1+q^{2n}) 
a_{i,n}z^{-n}],
$$
and 
$$k_i^-(z)= \exp[(q-q^{-1}) \Sigma_{n>0}-q^{2(n+c/2}/(1+q^{2n}) a_{i,n}
 z^{-n}].$$
Let 
$$\bar x^+(z)=x^+(z)k_i^-(z),$$ and 
$$\tilde x^+(z)=k_i^+(z) x^+(z).$$ 

\begin{thm}
$$\bar x_i^+(z)\bar x_i^+(w)=\bar x_i^+(w)\bar x_i^+(z), $$
$$\tilde x_i^+(z)\tilde x_i^+(w)=\tilde x_i^+(w)\tilde x_i^+(z). $$
\end{thm}

It is obvious that the set of current operators
$\vphi_i(z)$, $\psi_i(z)$, $\tilde  x_i^+(z)$ and 
$x^-(z)$ and the the set 
of the current operators 
$\vphi_i(z)$, $\psi_i(z)$, $\bar x_i^+(z)$ and 
$x^-(z)$  generate the quantum affine algebra 
$U_q(\hat {\frak sl}(n))$ respectively. 
The reformulation of  the quantum affine algebra 
$U_q(\hat {\frak sl}(n))$ with current operators 
$\vphi_i(z)$, $\psi_i(z)$, $\bar x_i^+(z)$ and 
$x^-(z)$  is the key for the quantized semi-infinite 
construction in the next section, namely we need to use the 
kernel coming from the current operator  $\bar x_i^+(z)$ 
to define the semi-infinite space. 
From now on, we will restrict ourselves to 
the case of $U_q(\hat {\frak sl}(2))$. The case for 
$U_q(\hat {\frak sl}(n))$ can be dealt with in a similar way
\cite{FS}. 

For the case of  $U_q(\hat {\frak sl}(2))$,
the  relations between 
$\psi(z)$, $\bar x^+(z)$ is the same as that of $\psi_i(z)$, $x^+(z)$, 
however the rest are changed, which we will write then  down below. 

\begin{prop}
\begin{align*}
& \vphi(z)\bar x ^{+}(w)\vphi(z)^{-1}=
 f_1(z/w) g(z/wq^{ -\frac{1}{2}c})^{+1}\bar x^{+}(w), \\
& (f_2(w/z)\bar x^+(z)x^-(w)-x^-(w) 
\bar  x^+(z)= \\
&\frac{\delta_{i,j}}{q-q^{-1}}
  \left\{ \delta(z/wq^{-c})\psi(wq^{\frac{1}{2}c})k^-(z)-
          \delta(z/wq^{c})\vphi(zq^{\frac{1}{2}c})k^-(z) \right\}, \\
&\bar x^+(z)\bar x^+(w)= \bar x^+(w)\bar x^+(z), 
\end{align*}
where $f_1(z/w)= (\frac{(1-z/wq^2q^{c/2})(1-z/wq^{c/2})}
{(1-z/wq^2q^{3c/2})(1-z/wq^{3c/2})})^{-1}$ and $f_2(w/z)= (q^2q^cw/z-1)(w/zq^c-
1). $
\end{prop}

Similarly, one can write down the relations between 
$\vphi(z)$, $\psi(z)$, $\tilde  x^+(z)$ and 
$x^-(z)$, which we will omit here.  In the next section, we will 
$\bar  x^+(z)$ and 
instead of  $ x^+(z)$
 as the the current operator  
for  our semi-infinite construction of representations of 
 $U_q(\hat {\frak sl}(2))$.

\section{Quantum integrability condition and 
semi-infinite construction}

The integrability condition of the current operator $e(z)$ induces the 
semi-infinite construction for the unquantized case. 
The quantum integrability condition was studied in 
\cite{DM}, which is stated as the following: 

\begin{thm}
For any level $k\geq1$ integrable module of
$U_q\left(\hat{\frak sl}(2)\right)$, 
the correlation functions of $x^+(z_1)x^+(z_2)...x^+(z_k)x^+(z_{k+1})$ 
is zero if $z_2/z_1=z_3/z_2=\ldots=z_{k+1}/z_k=q^2$,
the correlation functions of $x^-(z_1)x^-(z_2)...x^-(z_k)x^-(z_{k+1})$ 
is zero
if $z_1/z_2=z_2/z_3=\ldots=z_k/z_{k+1}=q^2$ 
\end{thm}

However this condition can not be directly used 
for the semi-infinite construction,  because the noncommutativity
of the current operator $x^+{z}$ of $U_q\left(\hat{\frak sl}(2)\right)$. 
However the theorem above  implies:

\begin{cor}
For any level $k\geq1$ integrable module of
$U_q\left(\hat{\frak sl}(2)\right)$, 
$$\bar x^+(z_1)\bar x^+(z_2)...\bar x^+(z_k)\bar x^+(z_{k+1})=0$$ 
if $z_2/z_1=z_3/z_2=\ldots=z_{k+1}/z_k=q^2$, and 
\end{cor}

{\bf Proof.}
Let $\bar F(z_1,...,z_n)$  be the correlation function of
 a vector $v$ in   any level $k\geq1$ integrable module of
$U_q\left(\hat{\frak sl}(2)\right)$ and $v^*$ in the dual space 
of this level $K$ module, 
$<v^*, 
\bar x^+(z_1)\bar x^+(z_2)...\bar x^+(z_k)\bar x^+(z_{k+1})v>$
. Then we 
have that 
$$<v^*, 
\bar x^+(z_1)\bar x^+(z_2)...\bar x^+(z_k)\bar x^+(z_{k+1})v>= $$
$$<v^*, \prod_{i,j} \frac {(z_i-z_jq^2)}{z_i-z_j}
 x^+(z_1) x^+(z_2)... x^+(z_k) x^+(z_{k+1} 
k^+(z_1) k^+(z_2)... k^+(z_k) k^+(z_{k+1}v >.$$
Because 
$<v^*, 
 x^+(z_1) x^+(z_2)... x^+(z_k) x^+(z_{k+1})v_1>$ for 
for a vector $v_1$ in   any level $k\geq1$ integrable module of
$U_q\left(\hat{\frak sl}(2)\right)$ and 
the function $\prod_{i,j} \frac {(z_i-z_jq^2)}{z_i-z_j}$ is not  zero 
 if $z_2/z_1=z_3/z_2=\ldots=z_{k+1}/z_k=q^2$, we have 
 $$\bar F(z_1,...,z_n)=0. $$ 
Because of the commutativity of the operator 
$\bar x^+(z)$, we have 
$\bar x^+(z)\bar x^+(z)$ is a well-defined operator, thus 
$$\bar x^+(z)\bar x^+(z)=0.  $$

With the preparation above, in this section, we will describe a quantized 
semi-infinite construction along the line of \cite{FS}. Their starting point 
for the case of $\hat {\frak sl}(2)$ is the integrability condition for 
level $k$ integrable modules, namely on any level $K$ module from the 
category of representations with highest weight is a sum of 
of irreducible integrable  representations, if and only if 
$e(z)^{k+1}$ is zero. 

similarly we can make the following claim: 
\begin{thm}
Any level $K$ module  of  $U_q\left(\hat{\frak sl}(2)\right)$  from the 
category of representations with highest weight is a sum of 
of irreducible integrable  representations, if and only if 
$\bar x^+(z)\bar x^{(zq^2)}.....x^(zq^{2k})$ is zero. 
\end{thm}

{\bf Proof}
The theorem above already gives the proof for half of the 
theorem.  The other half comes form the fact that 
if we quotient the relation $q=1$, the condition
$\bar x^+(z)\bar x^{(zq^2)}.....x^(zq^{2k})$ is zero
 simply degenerates into the condition that 
$e(z)^{k+1}$ is zero. Thus, it is integrable as a module of 
$ \hat{\frak sl}(2)$. From the theory of Lusztig, we know that 
all the integrable highest weight module must comes from the 
corresponding quantized module. Thus the module is also an 
integrable module when $q$ is generic. 

We will start our semi-infinite construction with the 
irreducible integrable module $V_{0,1}$ with the the highest weight vector  
$v_{0,1}$, such that the weight of the highest weight vector  is 
$0$ and the central element $c$ acts as $1$.

Let $\bar x^+(z)= \Sigma \bar x^+_i z^{-i}$ and  $U(\bar x)$ be the subalgebra
generated by  $\bar x^+_i$. 
We denote  $U(\bar x)^-$ the subalgebra  generated by 
$\bar x^+(n), n\geq 0$ and and $U(\bar x)^+$ the subalgebra 
 generated by $\bar x^+(n), n<0$
Let $W=  U(\bar x)v_{0,1}$, because $U(\bar x)^+v_{0,1}=0$, 
we have that $W$ is equivalent to    $U(\bar x)^+/I)v_{0,1}$, where 
$I$ is an ideal. 

\begin{lemma}
The ideal $I $ is generated by $S_k^{1}=\Sigma \bar x_i\bar {x_{k-i}}(q^{2i}+ 
q^{2k-2i}$, for $k<-1$. 
\end{lemma}

{\bf Proof.}
From the quantum integrability condition above, 
 we know that the elements $ S_k^{1}$ for $k<-1$ are inside 
the ideal I. We will denote the ideal generated by those elements by 
$I'$. The proof follows from that fact that 
$U(\bar x)^+/I)v_{0,1}$ has the same character as the case, when 
quotient the relation $q-1=0$. Thus $I=I'$.

\begin{defn}
 $ \bar V_{0,1}$ is  an vector space with the basis of infinite mononomials 
$M$ in the form of $x_{i_1}x_{i_2}..x_{i_n}....$, where $\{ i_1, i_2,....\}$ 
is an infinite sequence of indices such that, for some $n$, $i_n$ is odd 
and $i_{p+1}=i_p+2$, if $p>n$. Let V be a quotient space of $\bar V$, the 
quotient is given by the following relations: 

(1) $\bar x_i$ and $\bar x_j$  commutes, if $i\neq j$, 

(2) if an element $m\in \bar V$ contains a part 
$x_ux_{2N+1}x_{2N+3}x_{2N+5}...   $ and 
$u>2N-1$, then  $m=0$.  

(3) The operator $S_k=\Sigma_{a+b=k} x_ax_b(q^{2b}+q^{2a})$ acts 
on $\bar V$ and $S_kv=0$ for $v\in \bar V$. 
\end{defn}

We  define the action of $\bar x_i$ simply  by multiplication. 
The action of $a_i$ for $i>0$ is given by 
$$a_i\bar x_{i_1}\bar x_{i_2}....= $$
$$ [a_i, \bar x_{i_1}]\bar x_{i_2}...+ 
\bar x_{i_1}[a_i, \bar x_{i_2}]\bar x_{i_3}....+$$ 
$$
...+ \bar x_{i_1}\bar x_{i_2}...[a_i, \bar x_{i_n}]\bar x_{i_{n+1}}....$$
This is an finite expression. 
We define the action of $a_0$ by that: 
$$a_0( \bar x_{2N+1}\bar x_{2N_3}\bar x_{2N+5}....)= -2N 
 \bar x_{2N+1}\bar x_{2N_3}\bar x_{2N+5}....$$
  The action of $a_{-1}$ is defined as: 
$$a_{-1}\bar x_1\bar x_3\bar x_5.......= $$
$$(\bar x_0\bar x_3\bar x_5......+ \bar x_1\bar x_2\bar x_5.....+\bar x_1\bar x_3\bar x_4\bar x_7\bar x_9....+
\bar x_1\bar x_3\bar x_5\bar x_6\bar x_9........+ ....)=$$
$$(\bar x_0\bar x_3\bar x_5.....-\frac{(q^6+1)}{(q^4+q^{2})}\bar x_0\bar x_3\bar x_5......+ 
\frac {(q^6+1)}{(q^4+q^{2})}\frac {(q^{10}+q^4)}{(q^6+q^8)}\bar x_0\bar x_3\bar x_5......=$$
$$ (\bar x_0\bar x_3\bar x_5....) 1/(1+\frac{(q^6+1)}{(q^4+q^{2})}). $$
Thus it converges if $|\frac{(q^6+1)}{(q^4+q^{2})})|<1$.

We would like to define the action of $a_{-2}$ as the following: 
$$a_{-2}\bar x_1\bar x_3\bar x_5.......= $$
$$(q^2-q^{-2})(\bar x_{-1}\bar x_3\bar x_5......+ 
\bar x_1\bar x_1\bar x_5.....+\bar x_1\bar x_3\bar x_3\bar x_7\bar x_9....+
\bar x_1\bar x_3\bar x_5\bar x_5\bar x_9........+ ....)=$$
$$(q^2-q^{-2})(\bar x_{-1}\bar x_3\bar x_5.....-(\frac{(q^4+1)}{(2q)}\bar x_0\bar x_2\bar x_5......+
 \frac{(q^6+q^{-2})}{(2q)}\bar x_{-1}\bar x_3\bar x_5......)-$$ 
$$(\frac {(q^{10}+q^2)}{(q^6+q^6)}\bar x_1\bar x_1\bar x_5\bar x_7......+
\frac {(q^{8}+q^4)}{(q^6+q^6)} \bar x_1\bar x_2\bar x_4\bar x_7...$$
To use the relation (1) (2) (3) to reduce 
this expression to prove the convergence of the expression is 
very complicated. Similar problems appears for the defining  the 
action of $a_{-n}, n<-2$. 

Thus we will use the same trick  played in \cite{FS}. Let 
 $U(\bar x, r)$ be the subalgebra  generated by $\bar x_i$. 
Let $\bar V_{0,1}(r)$ be the subspace $\bar V_{0,1}$, which consists of 
the element  $\bar x_{i_1}......\bar x_{i_n}......$ and $i_j>r$.

\begin{lemma} 
$\bar V_{0,1}(r)$ spans the whole space $\bar V_{0,1}$. 
\end{lemma} 

{\bf Proof}
The proof is the same as lemma 2.5.1 in \cite{FS}. The way to prove it is to
use the relation (3) to express any element in $\bar V_{0,1}$ with 
linear expression of elements in $\bar V_{0,1}$.

For any element expressed in a linear combination of elements in 
$\bar V_{0,1}(r)$, we define the action of $x^-(k)$, for $k+r>0$, as
that of $a_{-1}$ by using the commutation relations between 
$\bar x^+(z)$ and $\bar x^-(z)$. Because  $k+r>0$,, we know that it is well 
defined. As in  \cite{FS}, this is a well defined action, namely if we 
express an  element in two different ways in 
 $\bar V_{0,1}(r)$, the actions of $x^-(k)$ defined above coincide. 
 Again, with the commutation relation between 
$\bar x^+(z)$ and $\bar x^-(z)$, we can define the action of $a_{n},n<-1$, 
because $\bar x^+(z)$ and $\bar x^-(z)$ generate the whole algebra. 
Thus we have

\begin{thm}
There exists an action of  $U_q\left(\hat{\frak sl}(2)\right)$  on the 
space $\bar V_{0,1}$, such that $\bar V_{0,1}$ is equivalent to 
$ V_{0,1}$ as a representation of  $U_q\left(\hat{\frak sl}(2)\right)$ and 
the action of $\bar x_i$ acts by multiplication. 
\end{thm}

Let $\bar W$ be the set of the elements  $\bar x_{i_1}......\bar x_{i_n}......$
 in  $\bar V_{0,1}$,  such that  
$i_{j+1}-i_{j}>1$. 

\begin{prop}
$\bar W$ forms an linear independent basis of the space $\bar V_{0,1}$.
\end{prop}

The proof is the same as in \cite{FS}, which gives the character of the 
representation. 

Similarly as in \cite{FS}, a functional model for the description of 
 $W^*$, the dual space of $W$,   can be derived from the lemma above. 

As a commutative algebra, 
 $U(\bar x)$ can be identified with the space $C[t,t^{-1}]$. 
Let $U(\bar x)^+= \Sigma \oplus U(\bar x)^+(n)$, where
$ U(\bar x)^+(n)$ consists of the elements 
 $\bar x^+_{i_1}\bar x^+_{i_2}...\bar x^+_{i_n}$.
We identify any element    $\bar x^+_{i_1}\bar x^+_{i_2}...\bar x^+_{i_n}$ 
in $ U(\bar x)^+(n)$ as $t_1^{i_1}....t_n^{i_n}$, 
 where  $t_i$ are variables.  
Similarly we can any element $\bar x^+_{i_1}\bar x^+_{i_2}...\bar x^+_{i_m}$ 
in $W$ as $t_1^{i_1}....t_m^{i_m}$. Let $S^n(\Omega^1\Bbb C)$ be the 
space of expressions $f(t_1,.....t_n)dt_1....dt_n$, such that 
$f(t_1,.....t_n)$ is a symmetric functions and different 
$dt_i$ commute. $S^n(\Omega^1\Bbb C)$ is also called the space of 
$n$ particles. We can pair $S^n(\Omega^1\Bbb C)$ with 
$ U(\bar x)^+(n)$ by the following: 
$$<f(t_1,.....t_n)dt_1....dt_n,  t_1^{i_1}....t_n^{i_n}>= 
\text{RES}_{t_1=..=t_n=0} (<f(t_1,.....t_n) t_1^{i_1}....t_n^{i_n}dt_1....dt_n.
$$
Thus $W^*= \Sigma \oplus W^*\bigcap  S^n(\Omega^1\Bbb C)$. 

\begin{thm}
$$
W^*\bigcap  S^n(\Omega^1\Bbb C)=\{ f(t_1,.....t_n)dt_1....dt_n:
f=0, \text{if} t_1=t_2q^2 \} .$$
\end{thm} 

Similarly we can present the semi-infinite constructions for 
the higher level cases.

Let $V_{l,k}$ be the  irreducible highest weight representation of 
$U_q\left(\hat{\frak sl}(2)\right)$ 
with the action of $c$ as $k$ and the highest weight is $l$ times the 
fundamental weight and $v_{0,1}$ be its highest weight vector. 
Let $W_{l,k}=  U(\bar x)v_{l,k}$, because $U(\bar x)^+v_{l,k}=0$, 
we have that $W_{l,k}$ is equivalent to    $U(\bar x)^+/I_{l,k})v_{l,k}$, 
where 
$I_{l,k}$ is an ideal. 

\begin{lemma}
The ideal $I_{l,k}$ is generated by $$S_i^{k+1}=
\Sigma_{\Sigma a_i=-i} x_{a_1}x_{a_2}...x_{a_{k+1}}
 (\Sigma_{\sigma \in s_{k+1}}(q^{\Sigma_{i=2,K+1}2{i-1}a_{\sigma i}}),$$
$i<-k$ and $\bar x_{-1}^{k-l+1}$.  
\end{lemma}

\begin{defn}
Let $\bar V_{l,k}$ be a space spaned by the elements of the 
expressions in the form of 
  $\bar x_{i_1}......\bar x_{i_n}\bar x_{2N}^l\bar x_{2N+1}^{k-l}
\bar x_{2N+2}^l\bar x_{2N+3}^{k-l}.....$, such that

(1) $\bar x_{i}$ commutes with $\bar x_{j}$,

(2) if an element $m\in \bar V$ contains a part 
$\bar x_{i}\bar x_{2N}^l\bar x_{2N+1}^{k-l}
\bar x_{2N+2}^l\bar x_{2N+3}^{k-l}.....$ , $i>2N-1$ or 
$\bar x_{i}\bar x_{2N+1}^{k-l}
\bar x_{2N+2}^l\bar x_{2N+3}^{k-l}.....$, $i>2N$, then  $m=0$.  

(3) The operator $S_k=\Sigma_{\Sigma a_i=n} x_{a_1}x_{a_2}...x_{a_{k+1}}
 (\Sigma_{\sigma \in s_{k+1}}(q^{\Sigma_{i=2,K+1}2{i-1}a_{\sigma i}})$ acts 
on $\bar V$ and $S_kv=0$ for $v\in \bar V$, where $S_{k+1}$ is the 
permutation group on $k+1$ numbers. 
\end{defn} 

\begin{thm}
On the space $\bar V_{l,k}$, there is an action of 
$U_q\left(\hat{\frak sl}(2)\right)$ , such that the action of 
$\bar x^+(n)$ is given by comultiplication. This representation is the 
irreducible highest weight representation of 
$U_q\left(\hat{\frak sl}(2)\right)$ 
with the action of $c$ as $k$ and the highest weight is $l$ times the 
fundamental weight. 
\end{thm}

Let $\bar W(l,k)$
 be the set of the elements  $\bar x_{i_1}......\bar x_{i_n}......$
 in  $\bar V_{l,k}$,  such that  
$i_{j+k}-i_{j}>1$.

\begin{prop}
$\bar W_(l,k)$ forms an linear independent basis of the space $\bar V_{l,k}$.
\end{prop}

As in the case of $V_{0,1}$, 
let $U(\bar x)^+= \Sigma \oplus U(\bar x)^+(n)$, where
$ U(\bar x)^+(n)$ consists of the elements 
 $\bar x^+_{i_1}\bar x^+_{i_2}...\bar x^+_{i_n}$.
We identify any element    $\bar x^+_{i_1}\bar x^+_{i_2}...\bar x^+_{i_n}$ 
in $ U(\bar x)^+(n)$ as $t_1^{i_1}....t_n^{i_n}$, 
 where  $t_i$ are variables.  
Similarly we can any element $\bar x^+_{i_1}\bar x^+_{i_2}...\bar x^+_{i_m}$ 
in $W_{l,k}$ as $t_1^{i_1}....t_m^{i_m}$. Let $S^n(\Omega^1\Bbb C)$ be the 
space of expressions $f(t_1,.....t_n)dt_1....dt_n$, such that 
$f(t_1,.....t_n)$ is a symmetric functions and different 
$dt_i$ commute. $S^n(\Omega^1\Bbb C)$ is also called the space of 
$n$ particles. We can pair $S^n(\Omega^1\Bbb C)$ with 
$ U(\bar x)^+(n)$ by the following: 
$$<f(t_1,.....t_n)dt_1....dt_n,  t_1^{i_1}....t_n^{i_n}>= 
\text{RES}_{t_1=..=t_n=0} (<f(t_1,.....t_n) t_1^{i_1}....t_n^{i_n}dt_1....dt_n.
$$
Thus $W^*= \Sigma \oplus W^*\bigcap  S^n(\Omega^1\Bbb C)$. 

\begin{thm}
$$
W^*\bigcap  S^n(\Omega^1\Bbb C)=\{ f(t_1,.....t_{n})dt_1....dt_{n}:
f=0, 
\text{if } t_1=t_2q^2....=t_{k+1}q^{2(k)} $$
$$\text{,  or }
t_1=...t_{k-l+1}=0 \} .$$
\end{thm}

In Section 2, we define the 
operator $\bar x_i^+(z)$, it is clear that this can be also applied to other 
cases. Our next step is to extend the semi-infinite to the 
cases of $U_q(\hat {\frak  g})$, where $g$ is a simple-laced simple 
Lie algebra, for which we need to define  the proper 
$x_\alpha^+(z)$ associated to the roots of ${\frak g}$. 
The simplest case  ${\frak g}={\frak sl}(3)$ will be the subject of a 
subsequent paper. 
On the other hand, this paper follows the algebraic theory developed 
in \cite{FS} \cite{FS1}. The semi-infinite constructions can be 
geometrically understood according to the structure of the 
corresponding infinite dimensional flag manifold and the infinite 
Shubert cells. The geometric interpretation of the quantized 
semi-infinite construction is still an open problem. This  is 
related to another 
immediate problem  to extend the explicit  construction of 
modular functor \cite{FS2} to the quantized case. 
which should leads us toward certain  quantization of conformal field theory. 
One more possible application of such a construction is to generaliza such a 
construction to more general cases. From the point view of the functional 
realization of the dual space, one natural generalization is to 
substitute the condition $x_1=x_2q^2$, which is generalization of the 
classical condition $x_1=x_2$,  with more general conditions, for 
example  $x_1=x_2q_1=x_3 q_3 $. One would like ask a questions 
that what kind of structures  those generalized spaces possiblely 
represent? We believe it is related with the recent work about 
generalization of the quantum affine algebras \cite{DI}, where 
those kind of new conditions should be satisfied for the 
quantum current operators. We hope our construction can help us 
to understand the structures of those new algebras in \cite{DI}, 
for which we have not been  able to give any  concrete  realization of the 
non-trivial integrable representations yet. 

\bigskip 

{\it Acknowledgment}

\medskip
\noindent
The authors thank T. Miwa for discussions.
J.D. is supported by the grant Reward research (A) 08740020 from the 
Ministry of Education of Japan.

\end{document}